\documentclass[12pt,twoside,reqno]{amsart}
\usepackage[colorlinks=false,citecolor=blue]{hyperref}
\usepackage{cleveref}
\usepackage{mathptmx, amsmath, amssymb, amsfonts, amsthm, mathptmx, enumerate, color,mathrsfs}
\usepackage{graphicx}

\usepackage{epstopdf}
\usepackage[T1]{fontenc}
\usepackage{mathptmx}

\newtheorem{theorem}{Theorem}[section]
\newtheorem{proposition}[theorem]{Proposition}
\newtheorem{lemma}[theorem]{Lemma}
\newtheorem{corollary}[theorem]{Corollary}
\theoremstyle{definition}
\newtheorem{remark}[theorem]{Remark}
\newtheorem{definition}[theorem]{Definition}
\newtheorem{example}[theorem]{Example}

\numberwithin{equation}{section}

\begin{document}
	
\title[Generalized derivations on a C$^{\ast}$-algebra]{On the equivalence of all notions of generalized derivations whose domain is a C$^{\ast}$-algebra}

\author[A. Hosseini]{Amin Hosseini}
\address[A. Hosseini]{Kashmar Higher Education Institute, Kashmar, Iran}
\email{a.hosseini@kashmar.ac.ir}

\author[A.M. Peralta]{Antonio M. Peralta}
\address[A.M. Peralta]{Instituto de Matem{\'a}ticas de la Universidad de Granada (IMAG), Departamento de An{\'a}lisis Matem{\'a}tico, Facultad de
	Ciencias, Universidad de Granada, 18071 Granada, Spain.}
\email{aperalta@ugr.es}

\author[S. Su]{Shanshan Su}
\address[S. Su]{%
	School of Mathematics, East China University of Science and Technology, Shanghai, 200237 China \\
	(Current address) Departamento de An{\'a}lisis Matem{\'a}tico, Facultad de
	Ciencias, Universidad de Granada, 18071 Granada, Spain.}
\email{lat875rina@gmail.com}


\begin{abstract} Let $\mathcal{M}$ be a Banach bimodule over an associative Banach algebra $\mathcal{A}$, and let $F: \mathcal{A}\to \mathcal{M}$ be a linear mapping. Three main uses of the term \emph{generalized derivation} are identified in the available literature, namely,  \begin{enumerate}[$(\checkmark)$]
\item $F$ is a generalized derivation of the first type if there exists a derivation $ d : \mathcal{A}\to \mathcal{M}^{**}$ satisfying $F(a b ) = F(a) b + a d(b),$ for all $a,b\in \mathcal{A}$.
\item $F$ is a generalized derivation of the second type if there exists an element $\xi\in \mathcal{M}^{**}$ satisfying $F(a b ) = F(a) b + a F(b) - a \xi b,$ for all $a,b\in \mathcal{A}$.
\item $F$ is a generalized derivation of the third type if there exist two (non-necessarily linear) mappings $G,H : \mathcal{A}\to \mathcal{M}$ satisfying $F(a b ) = G(a) b + a H(b),$ for all $a,b\in \mathcal{A}$.
\end{enumerate} These three types of maps are not, in general, equivalent. Although the first two notions are well studied when $\mathcal{A}$ is a C$^*$-algebra, their connections with the third one have not yet been explored. In this note we prove that every generalized derivation of the third type from a C$^*$-algebra $\mathcal{A}$ to a Banach $\mathcal{A}$-bimodule $\mathcal{M}$ is automatically continuous. We also show that every (continuous) generalized derivation of the third type from $\mathcal{A}$ to $\mathcal{M}$ is a generalized derivation of the first and second type. Consequently, the three notions coincide in this case. We also explore some concepts of generalized Jordan derivations on a C$^*$-algebra and establish some continuity properties for them. 
\end{abstract}

\maketitle

%------------------------------------------------------------------------------------%
\pagestyle{myheadings}
%
%\markboth{\centerline {}}
%{\centerline {}}
%\bigskip
%\bigskip
%------------------------------------------------------------------------------------%
%------------------------------------------------------------------------------------%

\section{Introduction}

Several notions of generalized derivations from an algebra $\mathcal{A}$ to a bimodule $\mathcal{M}$ have been considered in the literature since early nineties. Derivations are among the most studied maps in the literature. Recall that a linear mapping $d$ from an associative algebra $\mathcal{A}$ to an $\mathcal{A}$-bimodule  $\mathcal{M}$ is called a \emph{derivation} (respectively, a \emph{Jordan derivation}) if it satisfies Leibniz' rule, that is,  $$d(ab) = d(a)b + a d(b) \ \ \hbox{ (respectively, } d(a^2) = d(a) a + a d(a)\hbox{)}, \ \  \ (\forall a, b \in \mathcal{A}).$$ For example, if we fix an element $x_0\in \mathcal{M},$ the mapping $d_{x_0}:\mathcal{A} \rightarrow \mathcal{M}$ given by $d_{x_0} (a) = [a, x_0] = a x _0 - x_0 a$ is a derivation. Derivations which are expressed as finite sums of maps of the form $d_{x_0}$ are called \emph{inner derivations}. \smallskip

The available literature contains at least three different uses of the term ``\emph{generalized derivation}''. The first one, in chronological order, appears in a paper by Bre\v{s}ar published in 1991 (see \cite{Bresar91}). Keeping the notation above, a linear mapping $G :\mathcal{A}\to \mathcal{M}$ is a \emph{generalized derivation of the first type} if there exists a derivation $d :\mathcal{A}\to \mathcal{M}$ such that the identity \begin{equation}\label{eq generalized derivaiton 1} G (a b ) = G(a) b + a d(b) 
\end{equation} holds for all $a,b\in \mathcal{A}$. In the original definition, $G$ is only assumed to be additive. Along this paper, we shall say that a linear mapping $G$ from a Banach algebra $\mathcal{A}$ to a Banach $\mathcal{A}$-bimodule $\mathcal{M}$ is a \emph{generalized derivation of the first type} if there exists a derivation $d :\mathcal{A}\to \mathcal{M}^{**}$ such that the identity \eqref{eq generalized derivaiton 1} holds for all  $a,b\in \mathcal{A}$. Note that the identity in \eqref{eq generalized derivaiton 1} is equivalent to say that $G-d : \mathcal{A}\to \mathcal{M}^{**}$ is a \emph{left multiplier}, that is, $$(G-d)(a b) = G(a) b + a d(b) - d(a) b - a d(b) = (G-d)(a) b\ \ (a,b\in \mathcal{A}).$$ Right multipliers can be similarly defined. This first notion is the one employed, for example, in the study by Heller, Miller, Pysiak, and Sasin connecting differential geometry, generalized derivations, and general relativity \cite{Hel}. \smallskip

Keeping the chronological order, the second definition of a generalized derivation was introduced by Nakajima in 1999 (cf. \cite[$(1.3)$]{Nak1999}). A reformulation of this notion was considered by Leger and Luks in \cite{LeLuks2000}. This is the notion employed by Alaminos, Bre\v{s}ar, Extremera, and Villena in the study of bounded linear operators preserving zero products (see \cite[Definition 4.1]{AlBreExVill09}), and for example in \cite{AyuKudPe2014,EP18,LiPan,LiZhou2010} and \cite{Shesta2012}. A \emph{generalized derivation of the second type} from a Banach algebra $\mathcal{A}$ into a Banach $\mathcal{A}$-bimodule $\mathcal{M}$ is a linear mapping $G:\mathcal{A}\to \mathcal{M}$ for which there exists $\xi\in  \mathcal{M}^{**}$ satisfying \begin{equation}\label{eq gener der of second type} G(ab) = G(a)  b + a  G(b) - a \xi b \hbox{ ($a, b \in A$).}
\end{equation} 

Every derivation is a generalized derivation of the second type, though the class determined by the latter maps is strictly wider (e.g. for each non-zero $a\in \mathcal{A}$, the mapping $x\mapsto x\circ a= \frac12 (a x + x a)$ is a generalized derivation of the second type on $\mathcal{A}$  which is not a derivation). Let us observe that if $\mathcal{A}$ is unital, with unit $\mathbf{1}$, module products of the from  $a \xi \mathbf{1}$ and $\mathbf{1} \xi b$ lie in $\mathcal{M}$ for all $a,b\in \mathcal{A}$, and hence the left and right multiplication operators $L_{\mathbf{1} \xi},R_{\xi \mathbf{1}}: a\mapsto \mathbf{1} \xi a, a \xi \mathbf{1}$ define two bounded linear operators from $\mathcal{A}$ to $\mathcal{M}$ (i.e., $\mathbf{1} \xi$ and $\xi \mathbf{1}$ behave like a multiplier). In this case, the mapping $d= G-L_{\mathbf{1} \xi}:\mathcal{A}\to \mathcal{M}$ is a derivation and $G (a b ) = G(a) b + a d(b)$ for all $a,b\in \mathcal{A}$. Therefore, every generalized derivation of the second type is a generalized derivation of the first type. \smallskip

Suppose now that $\mathcal{A}$ is a unital algebra, and $G: \mathcal{A}\to \mathcal{M}$ is a generalized derivation of the first type with associated derivation $d$. Since $G-d$ is a left multiplier, we have $(G-d) (a) = (G-d) (\mathbf{1}) a $ for all $a\in \mathcal{A}$, and thus $$G(a b ) = G(a) b + a G(b) - a (G-d) (b) =  G(a) b + a G(b) - a (G-d) (\mathbf{1}) b  \ \  (a,b\in\mathcal{A}),$$ which shows that generalized derivations of the first and second type coincide in this case. \smallskip

Furthermore, every bounded left multiplier $L$ from a C$^*$-algebra $\mathcal{A}$ to a Banach $\mathcal{A}$-bimodule $\mathcal{M}$, is of the form $L(a) = \xi a,$ where $\xi\in M^{**}$ satisfies $\xi \mathcal{A}\subseteq \mathcal{A}$ (cf. \cite{AlBreExVill09,AyuKudPe2014}).\label{ref gd of type 1 are of type 2} Therefore, every continuous generalized derivation of the first type from $\mathcal{A}$ to $\mathcal{M}$ is a generalized derivation of the second type. Conversely, if $G: \mathcal{A}\to \mathcal{M}$ is a generalized derivation of the second type with respect to an element $\xi\in \mathcal{M}^{**}$, the mapping $d = G-L_{\xi} : \mathcal{A}\to \mathcal{M}^{**}$ is a derivation, and $G(a b ) = G(a) b + a d(b)$ for all $a,b\in \mathcal{A}$. We have thus shown that continuous generalized derivations of the first and second type from $\mathcal{A}$ to $\mathcal{M}$ coincide. It is important to note that continuity has been assumed to establish the equivalence between generalized derivations of the first two types from a C$^*$-algebra $\mathcal{A}$ to a Banach $\mathcal{A}$-bimodule. It should also be mentioned that every generalized derivation of the second type from a C$^*$-algebra into a Banach bimodule is automatically continuous (cf. \cite[Proposition 2.1]{JamPeSidd2015}).\smallskip

The third notion of generalized derivation is motivated by a definition introduced by Jimen{\'e}z-Gestal and P{\'e}rez-Izquierdo in \cite{J1}, and employed by Shestakov in \cite{Shesta2012, Shesta2014}. However, here we relax some linearity assumptions included in the original statement.   

\begin{definition} Let $\mathcal{A}$ be an algebra, and let $\mathcal{M}$ be an $\mathcal{A}$-bimodule. A ternary derivation from $\mathcal{A}$ to $\mathcal{M}$ is a $3$-tuple $(F,G,H)$, where $G, H:\mathcal{A} \rightarrow \mathcal{M}$ are two (non-necessarily linear) mappings, and $F:\mathcal{A} \rightarrow \mathcal{M}$ is a linear map satisfying $F(ab) = G(a)b + a H(b)$ for all $a, b \in \mathcal{A}$. We shall also say that $F:\mathcal{A} \rightarrow \mathcal{M}$ is a ternary derivation with associated mappings $G, H:\mathcal{A} \rightarrow \mathcal{M}$. Since the term ``ternary derivation'' is also employed in other settings with another meaning  (like in the case of JB$^*$-triples), we shall better say that the triplet $(F,G,H)$ behaves like a derivation. In this case, the mapping $F$ is called a \emph{generalized derivation of the third type}.  
\end{definition} 

There exist examples of $3$-tuples $(F,G,H): \mathcal{A}\to \mathcal{M}$ behaving like a derivation where $G$ and $H$ are not necessarily linear (see Example~\ref{example triple behaving as a triple derivation without linearity}).\smallskip

If $D$ is a derivation from an algebra $\mathcal{A}$ to an $\mathcal{A}$-bimodule $\mathcal{M}$, the triplet $(D,D,D)$ behaves like a derivation. Therefore every derivation is a generalized derivation of the third type.\smallskip

Komatsu and Nakajima presented in \cite{KoNak2003} a detailed study on the relationships among the different notions of generalized derivations of first, second, and third type, and their formal properties; but mainly in the unital case and from an algebraic perspective. It is perhaps worth recalling some basic connections. Suppose $F: \mathcal{A}\to \mathcal{M}$ is a generalized derivation of the second type satisfying \eqref{eq gener der of second type}. By assuming that $\mathcal{A}$ is unital, it follows that the element $\mathbf{1} \xi \mathbf{1} = \mathbf{1} G(\mathbf{1}) \mathbf{1}\in \mathcal{M}$. The expression 
$$ F(a b) = (F-R_{\mathbf{1} \xi \mathbf{1}}) (a) b + a F(b) = F (a) b + a (F-L_{\mathbf{1} \xi \mathbf{1}}) (b),$$ is valid for all $a,b\in \mathcal{A}$, and hence the triplets  $(F,F-R_{\mathbf{1} \xi \mathbf{1}},F)$ and $(F,F,F-L_{\mathbf{1} \xi \mathbf{1}})$ behave like derivations. That is, every generalized derivation of the second type is a generalized derivation of the third type. Conversely, it was already observed by Shestakov in \cite[Lemma 1]{Shesta2012} that if $\mathcal{A}$ is a unital associate algebra, every generalized derivation of the third type $F: \mathcal{A}\to \mathcal{A}$ is a generalized derivation of the second type (see also \cite[Lemma 4.1 and Corollary 4.5]{KoNak2003}). The argument works in our general setting too. Suppose $(F,G,H)$ is a $3$-tuple of mappings from a unital associative algebra $\mathcal{A}$ into a unital $\mathcal{A}$-bimodule, $\mathcal{M},$ behaving like a derivation, with $F$ being linear. As we shall see in Lemma~\ref{l unital case}, the mappings $G \mathbf{1}$ and $\mathbf{1} H$ are linear and satisfy $G(a)\mathbf{1} = F(a) - a \mathbf{1} H(\mathbf{1})$ and $\mathbf{1} H(a) = F(a) - G(\mathbf{1}) a,$ for all $a\in \mathcal{A}$. Therefore, it follows that
$$F (a b ) = F(a) b + a F(b) - a \Big( H(\mathbf{1}) +   G(\mathbf{1})\Big) b, \hbox{ for all } a,b\in \mathcal{A},$$ which shows that $F$ is a generalized derivation of the second type. \smallskip
 
However, if we relax the assumption that $\mathcal{A}$ is unital, it is not clear whether every generalized derivation of the third type is of the second type. This naturally gives rise to a question: are these three notions of generalized derivations from a (non-necessarily unital) C$^*$-algebra $\mathcal{A}$ into an $\mathcal{A}$-bimodule equivalent with each other? In this paper, we shall give an affirmative answer to this question. In fact, this problem is equivalent to a question on automatic continuity.\smallskip

Bounded linear operators which are generalized derivations admit a useful algebraic characterization. 

\begin{theorem}\label{t characterization bli generalized derivations}{\rm(\cite[Theorem 4.5]{AlBreExVill09}, \cite[Proposition 4.3]{BurFerPe2014}, \cite[Theorem 2.11]{AyuKudPe2014})} Let $T:\mathcal{A}\to \mathcal{M}$ be a bounded linear operator from a C$^*$-algebra to an essential Banach $\mathcal{A}$-bimodule. Then the following statements are equivalent:
	\begin{enumerate}[$(a)$]\item $T$ is a generalized derivation {\rm(}of the second type{\rm)}.
		\item $a T(b) c =0,$ whenever $ab =bc =0$ in $\mathcal{A}$.
		\item $aT(b)c =0,$ whenever $a b =b c =0$ in $\mathcal{A}_{sa}$.
	\end{enumerate}
\end{theorem} 

Recall that a Banach $\mathcal{A}$-bimodule $\mathcal{M}$ is called \emph{essential} if the linear span of the set $\{ a x b: a,b\in A, x\in \mathcal{M}\}$ is dense in $\mathcal{M}$. It is worth noting that in case that $\mathcal{A}$ is a Banach algebra admitting a bounded approximate unit $(u_j)_j$ and $\mathcal{M}$ is essential, we can deduce from the celebrated Cohen factorization theorem (see \cite[Theorem 32.22]{HewRosVolII} and \cite[Theorem in page 108, \S~1.15]{CLM1979}) that, for each $x\in \mathcal{M}$ the nets $(u_j x)_j$ and $(x u_j)_j$ converge in norm to $x$.\smallskip      

In case of a linear mapping acting on a von Neumann algebra $\mathcal{W}$ and behaving like a generalized derivation at certain points of the domain, continuity becomes an automatic property. More concretely, let $T:\mathcal{W}\to \mathcal{W}$ be a linear mapping on a von Neumann algebra. Suppose that for each $a, b, c$ in any commutative von Neumann subalgebra $\mathcal{B}\subseteq \mathcal{W}$ with $ab =bc =0$ we have $aT(b)c =0$. Then $T$ is automatically continuous \cite[Theorem 2.12 and Corollary 2.15]{EP18}. Furthermore, the above statements $(a)$--$(c)$ are also equivalent to the next:
	\begin{enumerate}[$(a)$] %\item $T$ is a generalized derivation (and hence continuous).
		%\item $a T(b) c =0$, whenever $ab =bc =0$ in $\mathcal{W}$.
		%\item $a T(b)c =0$, whenever $ab =bc =0$ in $\mathcal{W}_{sa}$.
		\item[$(d)$] $aT(b)c +cT(b)a =0$, whenever $ab =bc =0$ in $\mathcal{W}_{sa}$.
		\item[$(e)$] $a T(b) a =0,$ whenever $ab =0$ in $\mathcal{W}_{sa}.$	
	\end{enumerate}

Every von Neumann algebra is unital, and hence every generalized derivation of the third type on a von Neuman algebra is a generalized derivation of the first and second type.\smallskip

The list of studies exploring the automatic continuity of derivations and related operators is quite wide. The pioneering theorems by Sakai \cite{S} and Ringrose \cite{R} prove that every derivation on a C$^{\ast}$-algebra and every derivation from a C$^{\ast}$-algebra $\mathcal{A}$ into a Banach $\mathcal{A}$-bimodule is automatically continuous. Hou and Ming \cite{Hu} proved that if $\mathcal{X}$ is a simple Banach space, and $\sigma, \tau: B(\mathcal{X}) \rightarrow B(\mathcal{X})$ are surjective and continuous at $0$, then every $(\sigma, \tau)$-derivation from $B(\mathcal{X})$ into itself is continuous. Recall that, if $\mathcal{A}$ is an algebra and $\sigma, \tau :\mathcal{A} \rightarrow \mathcal{A}$ are two mappings, a \emph{$(\sigma, \tau)$-derivation on $\mathcal{A}$} is a linear mapping $d: \mathcal{A} \rightarrow \mathcal{A}$ satisfying $$d(ab) = d(a)\sigma(b) + \tau(a) d(b), \hbox{ for all } a, b \in \mathcal{A}.$$ More results on automatic continuity can be found in 
\cite{EP18,GP13,GH18,H4,H5,H6,J*,J87,PeRu2014,S1}.\smallskip

Let us explain how, in the setting of C$^*$-algebras and essential Banach bimodules, the problem of determining whether every generalized derivation of the third type is of the second type is a problem of automatic continuity. Suppose $F, G, H$ are three mappings from a C$^*$-algebra $\mathcal{A}$ to an essential $\mathcal{A}$-bimodule $\mathcal{M}$, such that the triplet $(F,G,H)$ behaves like a derivation, and let us assume that $F$ is continuous. For arbitrary $a,b,c\in \mathcal{A}_{sa}$ with $ a b = bc =0$, let us choose, via functional calculus, a decomposition of $b$ in the form $b = b^{+}- b^{-}$, with $b^{+}, b^{-}\geq  0$, $b^{+} b^{-} =0,$ $a b^{+} = a b^{-} =0$, $b^{+} c=0$, and  $b^{-} c =0.$ By taking $d = \left(b^{+}\right)^{\frac12} + i \left( b^{-}\right)^{\frac12}$, we have $ a d = d c =0$, $d^2 = b,$ and $$ a F(b) c =  a F(d^2) c = a G(d) d c + a d H(d) c = 0.$$ Thus, Theorem~\ref{t characterization bli generalized derivations} implies that $F$ is a generalized derivation of the second type. So it suffices to prove that every generalized derivation of the third type from a C$^*$-algebra into an essential Banach bimodule is continuous. In Proposition~\ref{p continuous gd 3type} we prove that every continuous generalized derivation of the third type from a C$^*$-algebra $\mathcal{A}$ to a (non-necessarily essential) Banach $\mathcal{A}$-bimodule is a generalized derivation of the first and second type. Theorem~\ref{t automatic cont of gener der of type3} completes the picture by showing that every generalized derivation of the third type from a C$^*$-algebra $\mathcal{A}$ to a Banach $\mathcal{A}$-bimodule $\mathcal{M}$ is continuous. Consequently, every generalized derivation of the third type from $\mathcal{A}$ to $\mathcal{M}$ is a generalized derivation of the second type (and of course, of the first type).\smallskip
  
Generalized Jordan derivations of the second type from a C$^*$-algebra $\mathcal{A}$ into a Banach $\mathcal{A}$-bimodule $\mathcal{M}$ have already been considered (see, for example, \cite{AyuKudPe2014,BurFerGarPe2014, BurFerPe2014}). A linear mapping $F :\mathcal{A} \rightarrow \mathcal{M}$ is a \emph{generalized Jordan derivation of the second type} if there exists an element $\xi\in \mathcal{M}^{**}$ satisfying  $$F(a \circ b) = F(a) \circ b + a \circ F(b) - U_{a,b} (\xi), \hbox{ for all } a, b \in \mathcal{A},$$ where $a \circ b =\frac12  (a b + b a)$ and $U_{a,b} (\xi) =\frac12 ( a \xi b + b \xi a)$. Here, we shall say that a linear mapping $F:\mathcal{A} \rightarrow \mathcal{M}$ is a \emph{generalized Jordan derivation of the third type} if there exist (non-necessarily linear) maps $G,H :\mathcal{A} \rightarrow \mathcal{M}$ such that $$ F( a\circ b) = G(a) \circ b + a \circ H(b), \hbox{ for all } a,b\in \mathcal{A}.$$ Theorem~\ref{t automatic cont generalize Jordan der of type 3} below shows that if $\mathcal{A}$ is a C$^{\ast}$-algebra, every generalized Jordan derivation of the third type from $\mathcal{A}$ into a Banach $\mathcal{A}$-bimodule $\mathcal{M}$ is continuous. Moreover, if we additionally assume that $\mathcal{M}$ is essential, then $F$ is a generalized derivation of the second type; and in certain cases (for example, when $\mathcal{M} = \mathcal{A}$ or $\mathcal{M}= \mathcal{A}^*$) generalized Jordan derivations of the third type will be generalized Jordan derivations of the second type whose associated element in $\mathcal{M}^{**}$ will commute with all elements in $\mathcal{A}$ (see Proposition~\ref{central elements generalized Jordan der}).

%%%%%%%%%%%%%%%%%%%%%%%%%%%%%%%
%\vspace{1cm}
\section{Automatic continuity of generalized derivations of the third type}

The main goal of this section is to establish a result assuring the automatic continuity of every generalized derivation of the third type from a C$^*$-algebra $\mathcal{A}$ into a Banach $\mathcal{A}$-bimodule.  We begin with an example of a $3$-tuple of maps behaving like a derivation in which at least one of the last two maps in the triplet is not linear. 

\begin{example}\label{example triple behaving as a triple derivation without linearity} Let $\mathcal{A}$ be an algebra, and let
\begin{align*}
\mathfrak{A} = \Bigg\{\left [\begin{array}{ccc}
a & b & c\\
0 & 0 & 0\\
0 & 0 & e
\end{array}\right ] \ : \ a, b, c, e \in \mathcal{A}\Bigg\}.
\end{align*}
Clearly, $\mathfrak{A}$ is an algebra with respect to the natural matrix product. Let $f:\mathcal{A} \rightarrow \mathcal{A}$ be a (non-necessarily linear) mapping. Define mappings $F, G, H:\mathfrak{A} \rightarrow \mathfrak{A}$ by $$F\Bigg(\left(\begin{array}{ccc}
a & b & c\\
0 & 0 & 0\\
0 & 0 & e
\end{array}\right)\Bigg) = \left(\begin{array}{ccc}
0 & a & b\\
0 & 0 & 0\\
0 & 0 & 0
\end{array}\right), \ G\Bigg(\left(\begin{array}{ccc}
a & b & c\\
0 & 0 & 0\\
0 & 0 & e
\end{array}\right)\Bigg) = \left(\begin{array}{ccc}
-a & f(b) & -c\\
0 & 0 & 0\\
0 & 0 & -e
\end{array}\right),$$ $$\hbox{ and } \ \  H\Bigg(\left [\begin{array}{ccc}
a & b & c\\
0 & 0 & 0\\
0 & 0 & e
\end{array}\right ]\Bigg) = \left [\begin{array}{ccc}
a & a+b & b+c\\
0 & 0 & 0\\
0 & 0 & e
\end{array}\right ].$$
Observe that $F$ and $H$ are linear while $G$ is non-necessarily linear.
A straightforward verification shows that $F(AB) = G(A)B + A H(B),$ for all $A, B \in \mathfrak{A},$
which means that $(F, G, H)$ behaves like a ternary derivation on $\mathfrak{A}$, and hence $F$ is a generalized derivation of the third type on $\mathcal{A}$.\smallskip

Let $D$ be a derivation on $\mathcal{A}$. Consider the Banach algebra $\mathcal{B} =\left\{ \left(\begin{array}{cc}
	a & 0 \\
	0 & 0 
\end{array}\right) : a\in \mathcal{A} \right\}$ and the $\mathcal{B}$-bimodule $\mathcal{M} = M_2 (\mathcal{A}) =\left\{ \left(\begin{array}{cc}
a & b \\
c & d 
\end{array}\right) : a,b,c,d\in \mathcal{A} \right\}$, with the natural operations. Define $F, G, H : \mathcal{B}\to M_2 (\mathcal{A})$, by $F  \left(\begin{array}{cc}
a & 0 \\
0 & 0 
\end{array}\right) := \left(\begin{array}{cc}
	D(a) & 0 \\
	0 & 0 
\end{array}\right),$ $G  \left(\begin{array}{cc}
a & 0 \\
0 & 0 
\end{array}\right) := \left(\begin{array}{cc}
D(a) & f(a) \\
0 & 0 
\end{array}\right),$ and $H  \left(\begin{array}{cc}
a & 0 \\
0 & 0 
\end{array}\right) := \left(\begin{array}{cc}
D(a) & 0 \\
f(a) & 0 
\end{array}\right).$ Clearly, $G$ and $H$ are non-necessarily linear, and it is easy to see that $(F,G,H)$ behaves like a derivation.
\end{example}

It follows from the previous counterexample that if the triplet $(F,G,H)$ behaves like a derivation, the mappings $G$ and $H$ need not be, in general, linear. However, if $\mathcal{A}$ is unital, we can replace them by linear mappings. 

\begin{lemma}\label{l unital case} Let $\mathcal{A}$ be an algebra, and let $\mathcal{M}$ be an $\mathcal{A}$-bimodule.  Suppose $(F, G, H)$ is a $3$-tuple of mappings from $\mathcal{A} $ to $\mathcal{M}$ behaving like a derivation. Assume that $\mathcal{A}$ is unital with unit $\mathbf{1}$. Then the mappings $G\cdot \mathbf{1}, \mathbf{1}\cdot H: \mathcal{A}\to \mathcal{M},$ $x\mapsto G(x) \mathbf{1}$, $x\mapsto \mathbf{1} H(x)$ are linear, and the $3$-tuple $(F,G\cdot \mathbf{1}, \mathbf{1}\cdot H)$ behaves like a derivation.  If $\mathcal{M}$ is unital, the mappings $G$ and $H$ are linear. 
\end{lemma}

\begin{proof} The conclusion is clear from the linearity of $F$, the identities $(G\cdot \mathbf{1}) (a) = F(a) - a H(\mathbf{1}),$ $(\mathbf{1}\cdot H)(a)  = F(a) - G(\mathbf{1}) a$ ($a\in \mathcal{A}$), and $$F(ab) = G(a) b + a H(b) = G(a) \mathbf{1} b + a \mathbf{1} H(b) = (G\cdot \mathbf{1}) (a)\ b + a\ (\mathbf{1} \cdot H)(b) \ (\forall a,b\in \mathcal{A}).$$
\end{proof}

Let us recall some well-known concepts. Let $\mathcal{A}$ be an algebra and let $\mathcal{M}$ be an $\mathcal{A}$-bimodule. The \emph{modular left annihilator of $\mathcal{A}$} is the set  $lann(\mathcal{A})_{\mathcal{M}}: = \{m_0 \in \mathcal{M} \ | \ m_0 \mathcal{A} = \{0\}\}.$ Similarly, the \emph{modular right annihilator of $\mathcal{A}$} defined by  $rann(\mathcal{A})_{\mathcal{M}}: = \{m_0 \in \mathcal{M} \ | \ \mathcal{A} m_0 = \{0\}\}.$ The \emph{modular annihilator of $\mathcal{A}$} is $ann(\mathcal{A})_{\mathcal{M}}: = rann(\mathcal{A})_{\mathcal{M}} \cap lann(\mathcal{A})_{\mathcal{M}}$. In particular, if $\mathcal{A} = \mathcal{M}$, we find the usual notions of left annihilator, right annihilator, and annihilator of the algebra $\mathcal{A},$ which are denoted by $lann(\mathcal{A})$, $rann(\mathcal{A}),$ and $ann(\mathcal{A})$, respectively. \smallskip

Along this section, the left and right multiplication operators by an element $a$ in an associative algebra $\mathcal{A}$ will be denoted by $L_a$ and $R_a,$ respectively. If $\mathcal{M}$ is an $\mathcal{A}$-bimodule, the mappings $L_a$ and $R_a$ will also stand for the corresponding left and right multiplication operators on $\mathcal{M},$ respectively.\smallskip

The next auxiliary lemma will be used in our arguments.    

\begin{lemma}\label{l Gb and aH} Let $\mathcal{A}$ be an algebra, and let $\mathcal{M}$ be an $\mathcal{A}$-bimodule.  Suppose $(F, G, H)$ is a $3$-tuple of mappings from $\mathcal{A} $ to $\mathcal{M}$ behaving like a derivation. Then the following statements hold:
\begin{enumerate}[$(i)$] 
\item The mappings $\Psi_a : \mathcal{A} \rightarrow \mathcal{M}$, $\Psi_a (b) = a H(b) =L_a H(b)$ and $\Gamma_a : \mathcal{A} \rightarrow \mathcal{M}$, $\Gamma_a (b) = G(b)a = R_a G(b)$ are linear for every $a \in \mathcal{A}$.
\item If $rann(\mathcal{A})_{\mathcal{M}} = \{0\}$ {\rm(}respectively, $lann(\mathcal{A})_{\mathcal{M}} = \{0\}${\rm)}, the mapping $G$ {\rm(}respectively, the mapping $H${\rm)} is linear.
\end{enumerate}

\end{lemma}

\begin{proof}$(i)$ We shall show that, for every $a \in \mathcal{A}$, the mapping $\Psi_a$ is linear. Given $a, b, c \in \mathcal{A}$ and $\lambda \in \mathbb{C}$, we have
\begin{align*}
F(a \lambda(b + c)) & = G(a)(\lambda b + \lambda c) + a H(\lambda (b + c)) \\ & = \lambda G(a) b + \lambda G(a)c + a H(\lambda (b + c)).
\end{align*}
On the other hand, since $F$ is a linear mapping, we have the following expressions:
\begin{align*}
F(a \lambda(b + c)) & = \lambda F(ab) + \lambda F(ac) \\
& = \lambda G(a) b + \lambda a H(b) + \lambda G(a) c + \lambda a H(c),
\end{align*}
for all $a, b, c \in \mathcal{A}$. Comparing the previous two equations we get that
\begin{align}\label{eq 2.1}
a H(\lambda(b + c)) -  \lambda aH(b) - \lambda aH(c) = 0,
\end{align}
which yields
$$ \Psi_a(\lambda(b+c)) = \lambda \Psi_a(b) + \lambda \Psi_a(c).$$
Hence, $\Psi_a$ is a linear mapping for any $a \in \mathcal{A}$. Similarly, one can show that the mapping $\Gamma_a : \mathcal{A} \rightarrow \mathcal{M},$ $\Gamma_a (t) = G(t)a$ is linear.\smallskip

$(ii)$ Assuming that $rann(\mathcal{A})_{\mathcal{M}} = \{0\}$, it follows from \eqref{eq 2.1} that $$ H(\lambda(b + c)) = \lambda H(b) + \lambda H(c),$$ for all $b, c \in \mathcal{A}$ and all $\lambda \in \mathbb{C}$, and thus $H(a + b ) = H(a) + H(b),$ for all $a, b \in \mathcal{A}$, $H(0) =0$ and $H(\lambda a) = \lambda H(a)$ for all $\lambda \in \mathbb{C}$. The other statement can be similarly obtained.
\end{proof}

%In the following remark, we recall some known conditions under which $lann(\mathcal{A})_{\mathcal{M}} $ $= \{0\} $ $= rann(\mathcal{A})_{\mathcal{M}}$.

\begin{example} Let $\mathbb{Z}$ be the set of all integers. Set $
		\mathcal{A} = \left\{\left [\begin{array}{cc}
			2n & 0\\
			0 & 2n
		\end{array}\right ] \ : \ n \in \mathbb{Z} \right\}.$
	It is clear that $\mathcal{A}$ is a non-unital ring. Let $
		\mathcal{M} = \left\{\left [\begin{array}{cc}
			i & j\\
			0 & k
		\end{array}\right ] \ : \ i, j, k \in \mathbb{Z} \right\}.$
	A straightforward verification shows that $lann(\mathcal{A})_{\mathcal{M}} = \{0\} = rann(\mathcal{A})_{\mathcal{M}}$.
\end{example}

If $\mathcal{A}$ is unital, $lann(\mathcal{A}) = \{0\} = rann(\mathcal{A})$. Also, if $\mathcal{A}$ is a semiprime algebra, then it is clear that $lann(\mathcal{A}) = rann(\mathcal{A}) = \{0\}$. If $\mathcal{A}$ satisfies \emph{Condition (P)} in \cite{F}, then it is routine to see that $lann(\mathcal{A}) = rann(\mathcal{A}) =\{0\}$; where an algebra $\mathcal{A}$ satisfies \emph{Condition (P)} if $a a_{0} a= \{0\}$ for any $a \in \mathcal{A}$ implies that $a_{0} = 0$.\smallskip

The dual space, $\mathcal{A}^*$, of a Banach algebra $\mathcal{A}$ is a Banach $\mathcal{A}$-bimodule with module operations given by $(\varphi a) (b) = \varphi (a b)$ and $(a \varphi)  (b) = \varphi (b a)$, for all $a,b\in \mathcal{A}$, $\varphi\in \mathcal{A}^*.$ If $\mathcal{A}$ satisfies a Cohen factorization type property (i.e. for every $c\in \mathcal{A}$ there exists $a,b\in \mathcal{A}$ with $c = a b $), we have $lann_{\mathcal{A}^*}(\mathcal{A})=\{0\} =rann_{\mathcal{A}^*}(\mathcal{A})$. Recall that, by Cohen's factorization theorem \cite[Corollary 2.26]{HewRosVolII}, every Banach algebra with a bounded left approximate unit satisfies such a  factorization property. To see the statement concerning modular annihilators, if $\varphi\in lann_{\mathcal{A}^*}(\mathcal{A}),$ we have $\varphi (a b) =0 $ for all $a,b\in \mathcal{A}$, and thus, the factorization property implies that $\varphi =0$. \smallskip

Our next proposition shows that every continuous generalized derivation of the third type from a C$^*$-algebra $\mathcal{A}$ into a Banach $\mathcal{A}$-bimodule is a generalized derivation of the first and second type. 

\begin{proposition}\label{p continuous gd 3type} Let $\mathcal{A}$ be a C$^*$-algebra, and let $F:\mathcal{A}\to \mathcal{M}$ be a continuous generalized derivation of the third type from $\mathcal{A}$ into a Banach $\mathcal{A}$-bimodule. Then $F$ is a generalized derivation of the first and second type.  
\end{proposition} 

\begin{proof} It is known that the product of $\mathcal{A}$ and the module products on $\mathcal{M}$ can be extended to a product on $\mathcal{A}^{**}$ and $\mathcal{A}^{**}$-bimodule operations on $\mathcal{M}^{**}$ via the first (or the second) Arens extensions \cite[Theorem 2.6.15$(iii)$]{DalesBook2000}, respectively. Since $\mathcal{A}$ is a C$^*$-algebra, its product is Arens regular, that is, the first and second Arens products on $\mathcal{A}^{**}$ coincide, and the latter is a von Neumann algebra with respect to this product (cf. \cite[Corollary 3.2.37]{DalesBook2000}). It is also known that the following properties holds: for each $a\in \mathcal{A}$, $\tilde{a}\in \mathcal{A}^{**}$, $x\in \mathcal{M}$ and $z\in \mathcal{M}^{**}$, the mappings $y\mapsto a y$, $y\mapsto  y \tilde{a}$ (respectively, $b\mapsto x b $, $b\mapsto  b z $) are weak$^*$ continuous maps on $\mathcal{M}^{**}$ (respectively, from $\mathcal{A}^{**}$ to $\mathcal{M}^{**}$) \cite[Proposition A.3.52]{DalesBook2000}; if $(a_\lambda)$ and $(x_\mu)$ are nets in $\mathcal{A}$ and $\mathcal{M}$, respectively, such that $a_\lambda \to a\in \mathcal{A}^{**}$ in the weak$^*$ topology of $\mathcal{A}^{**},$ and $x_\mu\to x\in \mathcal{M}^{**}$ in the weak$^*$ topology of $\mathcal{M}^{**}$, then \begin{equation}\label{eq product bidual module} a x =  \lim_{\lambda} \lim_{\mu} a_\lambda x_\mu \hbox{ and } x a =  \lim_{\mu} \lim_{\lambda}  x_\mu a_\lambda,
	\end{equation} in the weak$^*$ topology of $\mathcal{M}^{**}$ (cf. \cite[(2.6.26)]{DalesBook2000}).\smallskip

As observed in Lemma~\ref{l Gb and aH}, for each $a,b\in \mathcal{A}$, the maps $\Psi_a (\cdot) = a H(\cdot) = F(a \cdot) - L_{G(a)}$ and $\Gamma_b (\cdot) = G(\cdot) b = F(\cdot b ) - R_{H(b)}$ are linear, and in this case continuous by the assumptions on $F$. Let us consider the following sets of operators $$\Gamma = \left\{ G(\cdot) b : b\in \mathcal{A}, \ \|b\|\leq 1  \right\}, \hbox{ and } \Psi=\left\{ a H(\cdot)  : a\in \mathcal{A}, \ \|a\|\leq 1  \right\}.$$ For each $a\in \mathcal{A}$ we have $$\|   \left(G(\cdot) b\right) (a) \| = \|   G(a) b \| \leq \|F\| \|a\| \|b\| + \| a H(\cdot)\| \|b \| \leq \|F\| \|a\|+ \| a H(\cdot)\|,$$ and hence the uniform boundedness principle assures the existence of a positive $K_1$ satisfying $\|G(\cdot) b\| \leq K_1$ for all $b\in \mathcal{A}$ with $\|b\|\leq 1$. Similarly, there exists a positive $K_2$ satisfying $\|a H(\cdot) \| \leq K_2$ for all $a\in \mathcal{A}$ with $\|a\|\leq 1$.\smallskip

Let us take an approximate unit $(u_j)_j$ in $\mathcal{A}$. If we fix an element $a\in \mathcal{A},$ the net $(a u_j )_j$ converges in norm to $a$, and hence, by the continuity of $F$, $F(a u_j )$ tends to $F(a)$ in norm. It is also know that $(u_j)_j\to \mathbf{1}$ in the weak$^*$-topology of $\mathcal{A}^{**}$, where $\mathbf{1}$ stands for the unit in $\mathcal{A}^{**}$. We therefore deduce from \eqref{eq product bidual module} that $(G(a) u_j)_j \to G(a)\mathbf{1}$ in the weak$^*$-topology of $\mathcal{M}^{**}$. Thus, the identity $F (a u_j ) = G(a) u_j + a H(u_j)$ implies that 
\begin{equation}\label{eq a H(uj) converge} \left\{\begin{aligned}
&\hbox{ the net } (a H(u_j))_j \hbox{ converges, in the weak$^*$-topology of } \mathcal{M}^{**},  \\ 
&\hbox{ to some } R(a)\in \mathcal{M}^{**} \hbox{ and }  F(a) = G(a)\mathbf{1} + R(a).   
	\end{aligned}\right.
\end{equation} It is easy to check that the mapping $R : \mathcal{A}\to \mathcal{M}^{**},$ $a\mapsto R(a)$ is linear.  Moreover, since by the properties shown above $\|  a H(u_j) \| \leq K_2 \| a\|$ for all $j$, we obtain $\|R(a)\|\leq K_2 \|a\|$ for all $a\in \mathcal{A}$, and thus $R$ is continuous. \smallskip

Let us see another interesting property of the operator $R$. By definition and \eqref{eq product bidual module} we get 
$$R(a b ) = w^*\hbox{-}\lim_{j} ab H(u_j) =  a w^*\hbox{-}\lim_{j} b H(u_j) = a R(b) \ \ (\forall a,b\in \mathcal{A}),$$ which guarantees that $R$ is a right multiplier from $\mathcal{A}$ to $\mathcal{M}^{**}$.  Furthermore, by what we proved above, $\| R (u_j) \|\leq K_2$ for all $j$, and hence, by the weak$^*$-compactness of the closed unit ball of $\mathcal{M}^{**}$, there exists a subnet, denoted again by $(u_j)_j$, such that $\lim_{j} R(u_j) = \xi\in \mathcal{M}^{**}$ in the weak$^*$-topology of $\mathcal{M}^{**}$. Since $ R(a u_j) = a R(u_{j})$ for all $j$, and $(a u_j)_j \to a$ in norm, the continuity properties of the module operations on $\mathcal{M}^{**}$ and of the mapping $R$ give $R(a) = a \xi$ for all $a\in \mathcal{A}$.\smallskip 

Similar arguments show the existence of a left multiplier $L:\mathcal{A}\to \mathcal{M}^{**}$ and $\eta\in \mathcal{M}^{**}$ satisfying $L(a) = \eta a$ and $F(a) = L(a) + \mathbf{1} H(a)$ for all $a\in \mathcal{A}$. \smallskip

Finally, we have $$\begin{aligned}
F(a b ) &= G(a) b + a H(b) = G(a) \mathbf{1} b + a\mathbf{1} H(b) \\
&= (F(a)- R(a)) b + a (F(b)-L(b)) = F(a) b + a F(b) - a (\xi + \eta ) b,  
\end{aligned} $$ for all $a,b\in \mathcal{A}$. 
\end{proof}

We shall next state two classical arguments in results on automatic continuity which are considered here from a more general point of view.

\begin{lemma}\label{l continuity on a closed finite-codimensional subspace} Let $Z$ be a closed subspace of a Banach space $X$ such that $X/Z$ is finite-dimensional, and let $Y$ be a normed space. Suppose $F: X\to Y$ is a linear mapping whose restriction to  $Z$ is continuous. Then $F$ is continuous. 
\end{lemma}

\begin{proof} Since the quotient $X/Z$ is finite-dimensional, the subspace $Z$ is topologically complemented in $X$, that is, there exists a continuous linear projection $P: X\to X$ whose image is $Z$ and $Z^{\prime}=(Id-P) (X)$ is finite-dimensional (cf. \cite[Lemma 4.21]{RudinBookFA1991}). Clearly, by the finite-dimensionality of $Z^{\prime}$,  $F|_{Z^{\prime}}: Z^{\prime}\to Y$ is continuous. Since $F|_{Z}$ is continuous by hypothesis, and $F (x)= F (P(x)) + F ((Id-P)(x))=F|_{Z} (P(x)) + F|_{Z^{\prime}}((Id-P)(x)),$ for all $x \in X$, the conclusion is clear. 
\end{proof}

The next lemma is a consequence of the uniform boundedness principle.

\begin{lemma}\label{l separate cont} Let $F:\mathcal{A}\to X$ be a linear mapping from a C$^*$-algebra into a normed space. Suppose that for each $a\in \mathcal{A}$ the mappings $FL_a, F R_a : \mathcal{A}\to X$, $x\mapsto F(a x)$ and $x\mapsto F(x a)$ are continuous. Then $F$ is continuous.
\end{lemma}

\begin{proof} The desired conclusion is clear when $\mathcal{A}$ is unital, since in that case $F(a) = F L_{\mathbf{1}} (a)$. In the general case, we observe that the bilinear mapping $(a,b)\mapsto V(a,b):=F(a b)$ is separately continuous by hypothesis, so by the uniform boundedness principle, $V$ is jointly continuous. Thus, there exists a positive $K$ such that $\| F (a b)\| \leq K \ \| a\| \ \|b\|$ for all $a,b\in \mathcal{A}$. For each positive element $c\in \mathcal{A}$  with $\|c\|\leq 1$ there exists a positive $d$ with $d^2 = c$, $\|d\|^2 = \|d^2 \| = \|c\|\leq 1$. Therefore $\|F(c)\| =\|F(d^2)\| \leq K \|d\|^2 = K \|c\| \leq K$. In particular, $F$ is bounded on the closed unit ball of $\mathcal{A}$. 
\end{proof}

A cornerstone result in the theory of C$^*$-algebras, obtained by Cuntz in \cite{Cuntz1976}, asserts that a semi-norm $p$ on a C$^*$-algebra $\mathcal{A}$ which is bounded on each commutative self-adjoint subalgebra of $\mathcal{A}$, is bounded on the whole of $\mathcal{A}$. Cuntz' theorem has been employed in results on automatic continuity of derivations, for example, in \cite{PeRu2014} Russo and the second author of this note apply it to prove that every Jordan derivation from a C$^*$-algebra $\mathcal{A}$ to a Banach $\mathcal{A}$-bimodule is continuous and an associative derivation.  More recently, An and He employ Cuntz' theorem to prove that if $n\neq m$, then zero is the only $(m,n)$-Jordan derivation from a C$^*$-algebra into a Banach bimodule (cf. \cite{AnHe2019}). Next, we give another application by adapting the arguments in the proof of Theorem~\ref{t automatic cont of gener der of type3}.\smallskip

We are now ready to prove one of the main results of this section, which confirms that every generalized derivation of the third type from a C$^*$-algebra $\mathcal{A}$ into a Banach $\mathcal{A}$-bimodule is continuous.

\begin{theorem}\label{t automatic cont of gener der of type3} Let $\mathcal{A}$ be a  C$^{\ast}$-algebra, and let $\mathcal{M}$ be a Banach $\mathcal{A}$-bimodule. Suppose that $F: \mathcal{A} \to \mathcal{M}$ is a generalized derivation of the third type. Then $F$ is continuous. Consequently, every generalized derivation of the third type from $\mathcal{A}$ to $\mathcal{M}$ is a generalized derivation of the second type. 
\end{theorem}

\begin{proof} By the previously commented  theorem of Cuntz (see \cite[Theorem 1.1 or Corollary 1.2]{Cuntz1976}), there is no loss of generality in assuming that $\mathcal{A}$ is commutative.\smallskip
	
The proof will be presented in several steps.  Let us begin by setting $$I = \{a \in \mathcal{A} \ : \ F L_a \hbox{ is continuous}\} \hbox{ and } J = \{a \in \mathcal{A} \ : \ L_a H \hbox{ is continuous}\}.$$

The identity $F(a b) = a H(b) + G(a) b\ $  ($a,b\in \mathcal{A}$) implies that $F L_a (\cdot)$ is continuous if and only if $L_a H (\cdot)$ is a bounded linear operator, and thus $I = J$. Since $\mathcal{A}$ is commutative we also have $$I = \{a \in \mathcal{A} \ : \ F L_a = F R_a \hbox{ is continuous}\}.$$

We shall next show that $I$ is a closed ideal of $\mathcal{A}$. Namely, take $a\in I$ and $b\in \mathcal{A}$. Clearly, the mapping $F L_a  L_b:\mathcal{A} \rightarrow \mathcal{M}$,  $c\mapsto  F(a b c)$ is continuous, and so $ab = ba \in I$. This means that $I$ is an ideal of $\mathcal{A}$ (recall that $\mathcal{A}$ is commutative). It is well-known that $I$ must be self-adjoint \cite[Corollary 4.2.10]{K-R}. \smallskip 

Next we show that $I$ is norm-closed. Suppose that $(a_n)_n $ is a sequence in $I$ converging to an element $a\in \mathcal{A}$. By the equality $I = J$, in order to show that $a\in I$, it suffices to show that the mapping $L_a H:\mathcal{A} \rightarrow \mathcal{M}$ is continuous. It follows from Lemma~\ref{l Gb and aH}$(i)$ that the mapping $L_{c} H = :\mathcal{A} \rightarrow \mathcal{M}$ is linear for every $c \in \mathcal{A}$. Since $(a_n)$ is a sequence in $I$, the linear mapping $L_{a_{n}} H:\mathcal{A} \rightarrow \mathcal{M}$ is continuous for all $n \in \mathbb{N}$. It is clear that $\lim_{n \rightarrow \infty} L_{a_n} H(c) = \lim_{n \rightarrow \infty} {a_n}H(c) = a H(c) = L_a H(c),$ for every $c \in \mathcal{A}$, and hence, by the uniform boundedness principle, we obtain that $L_a H$ is a continuous linear mapping, and so $a \in J = I$.  \smallskip

It follows from the above arguments that the restricted mapping $F|_{I} : I \to \mathcal{M}$ satisfies the following property: for each $a\in I$, the mapping $F|_{I} L_{a} = F|_{I} R_{a} : I \to \mathcal{M}$, $x\mapsto F (a x) = F (x a)$ is continuous. Lemma~\ref{l separate cont} assures that $F|_{I} : {I}\to  \mathcal{M}$ is continuous.\smallskip

We shall next show that \begin{equation}\label{Step 4} \hbox{${\mathcal{A}}/{I}$ is a finite-dimensional C$^*$-algebra.}
\end{equation} 

Suppose, on the contrary, that ${\mathcal{A}}/{I}$ is an infinite-dimensional C$^{\ast}$-algebra. Another classical result in C$^*$-algebra theory (see \cite[Exercise 4.6.13]{K-R}) proves that in such a case there exists an infinite sequence $(b_n +I)_n$ of mutually orthogonal non-zero positive elements in ${\mathcal{A}}/{I}$, that is, $(b_n + I) (b_k+ I) = 0,$ for all $n \neq k$, $b_n\geq 0$ and $b_n\neq 0$ for all $n$. By \cite[Exercise 4.6.20]{K-R} we can always lift the sequence $(b_n +I)_n$ to a sequence $(c_n)_n$ of mutually orthogonal non-zero positive elements in  $\mathcal{A}$ satisfying $\pi(c_n) = b_n+I$, for all $n$, where $\pi:\mathcal{A} \rightarrow {\mathcal{A}}/{I}$ is the canonical projection. It follows from the fact that $\pi$ is a $^*$-homomorphism that $\pi(c_n^2)= \pi(c_n)^2 = (b_n+I)^2 \neq 0$ (observe that $\|  (b_n+I)^2 \| = \| b_n+I\|^2$), and thus $c_n^2 \notin I$ for all $n$. Therefore the mapping $F L_{{c^2_n}} : \mathcal{A} \rightarrow \mathcal{M}$ must be unbounded. By replacing $c_n$ with $\frac{c_n}{\|c_n\|}$, we can assume that $\|c_n\| =1$ for all $n$. By the unboundedness of the mapping $F L_{{c^2_n}}$, there exists $d_n\in \mathcal{A}$ satisfying $\|d_n\| \leq 2^{-n},$ and $\| F ({{c^2_n}} d_n)\| = \| F L_{{c^2_n}} (d_n)\| > \|G(c_n)\| + n$ for all natural $n$ --observe that the sequence $(\|G(c_n)\|)_n$ does not produce any obstacle here--. The elements $c_n$ and $d_n$ have been chosen to guarantee that the series $\sum_{n \geq 1} c_n d_n$ is (absolutely) convergent, and its limit $a_0 = \sum_{n = 1}^{\infty} c_n d_n$ satisfies $\|a_0\|\leq 1$ and $c_{m} a_0 = c_m^2 d_m$ for all $m\in \mathbb{N}$. By the hypotheses on $(F,G,H)$ we deduce that 
\begin{align*}
	\infty>\|H(a_0)\| &\geq \|c_m H(a_0)\|  = \|F(c_m a_0) - G(c_m) a_0\|  \geq \|F(c_m a_0)\| - \|G(c_m) a_0\| \\ & = \|F(c_m^2 d_m)\| - \|G(c_m) a_0\| \geq m + \|G(c_m)\| - \|G(c_m)\| = m,
\end{align*} for all natural $m$, which is impossible.\smallskip

Since $F|_{I}$ is continuous and $\mathcal{A}/I$ is finite-dimensional, we deduce from Lemma~\ref{l continuity on a closed finite-codimensional subspace} that $F$ is continuous. \smallskip

The final statement of the theorem is a consequence of Proposition~\ref{p continuous gd 3type}.
\end{proof}

The previous theorem generalizes the classical results by Ringrose \cite[Theorem 2]{R} and Sakai \cite{S}.  Let $F: \mathcal{A}\to \mathcal{M}$ be a linear mapping from a C$^*$-algebra to a Banach $\mathcal{A}$-bimodule. Suppose $F$ is a generalized derivation of the first type, that is, we can find a derivation $d: \mathcal{A}\to \mathcal{M}^{**}$ such that $F(a b ) = F(a) b + a d(b)\ \ $ ($a,b \in \mathcal{A}$). Clearly, the mapping $F:\mathcal{A}\to \mathcal{M}\subseteq \mathcal{M}^{**}$ is a generalized derivation of the third type, and thus continuous by \Cref{t automatic cont of gener der of type3}. Consequently every generalized derivation of the first type from a C$^*$-algebra to a Banach $\mathcal{A}$-bimodule is continuous.\smallskip

The following corollary summarizes some of the conclusions obtained up to now. 

\begin{corollary}\label{c all type of gen derivations coincide for Cstaralgebras} Let $F: \mathcal{A}\to \mathcal{M}$ be a linear mapping from a C$^*$-algebra to a Banach $\mathcal{A}$-bimodule. Then the following statements are equivalent:
\begin{enumerate}[$(a)$] 
	\item  $F$ is a generalized derivation of the first type, i.e., there is a derivation $d: \mathcal{A}\to \mathcal{M}^{**}$ such that $F(a b ) = F(a) b + a d(b)\ \ $ ($a,b \in \mathcal{A}$).
	\item $F$ is a generalized derivation of the second type, i.e., there is an element $\xi \in \mathcal{M}^{**}$ such that $F(a b ) = F(a) b + a F(b) - a \xi b\ \ $ ($a,b \in \mathcal{A}$).
	\item $F$ is a generalized derivation of the third type, i.e., there are two (non-necessarily linear) maps $G,H: \mathcal{A}\to \mathcal{M}$ such that $F(a b ) = G(a) b + a H(b)\ \ $ ($a,b \in \mathcal{A}$).
\end{enumerate} Furthermore, if any of the equivalent statements holds, the mapping $F$ is automatically continuous. 
\end{corollary}

Let ${Y}$ and ${Z}$ be Banach spaces and let $T: {Y} \rightarrow {Z}$ be a linear mapping. The \emph{separating space} of $T$ is the set $$\mathfrak{S}(T) = \Big\{z \in {Z} \ : \ \exists \ (y_n)_n \subseteq {Y} \ such \ that \  (y_n)_n \rightarrow 0, \hbox{ and } (T(y_n))_n \rightarrow z\Big\}.$$ By the closed graph theorem, $T$ is continuous if and only if $\mathfrak{S}(T) = \{0\}$. For additional information about separating spaces, the reader is referred to \cite{DalesBook2000}.\smallskip

We establish next some consequences of Theorem~\ref{t automatic cont of gener der of type3} and deduce some properties of the mappings appearing in a triplet behaving like a derivation.

\begin{corollary}\label{c 5}
Let $\mathcal{A}$ be a C$^{\ast}$-algebra, and let $\mathcal{M}$ be a Banach $\mathcal{A}$-bimodule. Suppose that $F, G, H:\mathcal{A} \rightarrow \mathcal{M}$ are three mappings, with $F$ linear, such that $(F,G,H)$ behaves like a derivation. The the following statements hold:	
\begin{enumerate}[$(i)$]\item  If $rann(\mathcal{A})_{\mathcal{M}} = \{0\} = lann(\mathcal{A})_{\mathcal{M}}$, the maps $F, G$ and $H$ are linear and continuous.
\item  If $\mathcal{M} = \mathcal{A}$, then $F, G$ and $H$ are bounded linear maps.
\item  If $\mathcal{M} = \mathcal{A}^*$, then $F, G$ and $H$ are bounded linear maps.
\end{enumerate}
\end{corollary}

\begin{proof}$(i)$ According to Theorem~\ref{t automatic cont of gener der of type3}, $F$ is a continuous generalized derivation of the second type, and it follows from Lemma~\ref{l Gb and aH}$(ii)$ that the mappings $G$ and $H$ are linear.  We shall show that $G$ and $H$ are continuous. Let $m_0 \in \mathfrak{S}(H) \subseteq \mathcal{M}$. Then, there exists a sequence $\{b_n\} \subseteq \mathcal{A}$ such that $\lim_{n \rightarrow \infty} b_n = 0$ and $\lim_{n \rightarrow \infty} H(b_n) = m_0$. For any arbitrary element $a \in \mathcal{A}$, we have
\begin{align*}
0 = \lim_{n \rightarrow \infty}F(a b_n) = \lim_{n \rightarrow \infty}(G(a) b_n + a H(b_n)) = a m_0,
\end{align*}
which means that $m_0 \in rann(\mathcal{A})_{\mathcal{M}}$. By hypothesis $m_0 = 0,$ and this implies that $H$ is continuous. Similarly arguments show that $G$ is continuous.\smallskip

$(ii)$ and $(iii)$ are straightforward consequences of $(i)$.
\end{proof}

Let $\mathcal{A}$ be an algebra and let $\mathcal{M}$ be an $\mathcal{A}$-bimodule. We recall that a linear mapping $F:\mathcal{A} \rightarrow \mathcal{M}$ is called an \emph{l-generalized derivation} if there exists a mapping $F_{l}: \mathcal{A} \rightarrow \mathcal{M}$ such that $F(ab) = F(a) b + a F_l(b)$ for all $a, b \in \mathcal{A}$. Similarly, $F$ is called an \emph{r-generalized derivation} if there exists a mapping $F_r: \mathcal{A} \rightarrow \mathcal{M}$ such that $F(ab) = F_r(a) b + a F(b)$ for all $a, b \in \mathcal{A}$. In \cite[Example 2.6]{H3}, we can find an example of an $r$-generalized derivation that is not an $l$-generalized derivation. Clearly, $l$- and $r$-generalized derivations are generalized derivations of the third type. The next result is therefore a straightforward consequence of Theorem~\ref{t automatic cont of gener der of type3} and Corollary~\ref{c 5}.

\begin{corollary}\label{c l-gen der} Let $\mathcal{A}$ be a C$^{\ast}$-algebra and let $\mathcal{M}$ be a Banach $\mathcal{A}$-bimodule. Then the following statements hold: \begin{enumerate}[$(i)$] 
\item  Every {l}-generalized derivation (respectively, {r}-generalized derivation) $F:\mathcal{A} \rightarrow \mathcal{M}$ is continuous.
\item Suppose that $rann(\mathcal{A})_{\mathcal{M}} = \{0\} = lann(\mathcal{A})_{\mathcal{M}}$. If $F:\mathcal{A} \rightarrow \mathcal{M}$ is an $l$- or $r$-generalized derivation with an associated mapping $d:\mathcal{A} \rightarrow \mathcal{M}$, then $F$ and $d$ are continuous linear mappings.
\end{enumerate}
\end{corollary}

In \cite{FreiMath1992} Freiberger and Mathieu considered an additional notion of ``\emph{generalized derivation}''.  Let $\mathcal{A}$ be a unital C$^*$-subalgebra of the unital C$^*$-algebra $\mathcal{B}$. A linear map $\delta : \mathcal{A} \to \mathcal{B}$ is called a ternary-generalized derivation if, for all $x, y, z \in \mathcal{A}$ we have $$\delta (xyz) = \delta(xy) z - x \delta(y) z + x \delta(yz).$$ These authors showed that $\delta$ is a ternary-generalized derivation if, and only if, there is a derivation $d:\mathcal{A} \to \mathcal{B}$ such that $\delta (xz) = d(x) z + x \delta(z)$ for all $x,z \in \mathcal{A}$. Obviously, ternary-generalized derivations are generalized derivations of the second and third type, and so our results can be also applied in this case.

%%%%%%%%%%%%%%%%%%%%%%%%%
%\vspace{1cm}
\section{Automatic continuity of generalized Jordan derivation of the third type}

Appropriate types of generalized Jordan derivations whose domain is a C$^*$-algebra are studied in this section. We shall mainly focus on the automatic continuity and the relationships between the different types. Contrary to the conclusions in the previous section, the natural definition for generalized Jordan derivations of the third type defines a strict subclass of the set of generalized Jordan derivations of the second type (cf. Proposition~\ref{central elements generalized Jordan der} and Remark~\ref{third type is stictly included in second}).  Generalized Jordan derivations of the second type have been considered in \cite{BurFerGarPe2014,BurFerPe2014,JamPeSidd2015}. Recall that a linear mapping $F$ from a Banach algebra $\mathcal{A}$ to a Banach $\mathcal{A}$-bimodule $\mathcal{M}$ is a \emph{generalized Jordan derivation of the second type} if there exists $\xi \in \mathcal{M}^{**}$ such that 
$$ F( a\circ b) = F(a) \circ b + a \circ F(b) - U_{a,b}(\xi), \hbox{ for all } a,b\in \mathcal{A},$$ where $a\circ b =\frac12 (a b + ba)$ and $U_{a,b} (\xi) = \frac12(a\xi b + b \xi a)$. If $\xi =0,$ we find the usual notion of Jordan derivation. Classic results by Ringrose and Johnson assure that every Jordan derivation from a C$^*$-algebra to a Banach bimodule is a derivation.\smallskip

There are some other attempts to define generalized Jordan derivations in the literature, According to \cite{BresBuk90}, an additive (in this note we shall assume linearity) mapping $T$ from a Banach algebra $\mathcal{A}$ into a left $A$-module $M$ is called a \emph{Jordan left derivation} if $T(a^2) = 2 a T(a)$ for every $a \in \mathcal{A}$. It is shown in \cite{BresBuk90} that the existence of non-zero Jordan left derivations from a prime ring $\mathcal{R}$ into a 6-torsion free left $\mathcal{R}$-module implies that $\mathcal{R}$ is a commutative ring. Vukman introduced in \cite{Vuk2008} the notion of $(m,n)$-Jordan derivation. %Let $m, n$ be two non-negative integers with $m+n \neq 0$. An additive (along this note linear) mapping $T$ from $\mathcal{A}$ into $\mathcal{M}$ is called an $(m, n)$-Jordan derivation if the identity $$ (m+ n) T(a^2) = 2 m a T(a) + 2 n T(a) a$$ holds for all $a\in \mathcal{A}$.
\smallskip

Jing and Lu coined for the first time the term generalized Jordan derivations in \cite{JingLu2003}. We shall say that a linear mapping $F$ from a Banach algebra $\mathcal{A}$ to a Banach $\mathcal{A}$-bimodule $\mathcal{M}$ is a \emph{generalized Jordan derivation} in the sense of Jing and Lu if there exists a (linear) Jordan derivation $\tau : \mathcal{A} \to \mathcal{M}$ satisfying $$F(a^2) = F(a) a + a \tau(a), \hbox{  for all } a \in \mathcal{A}.$$ By polarizing we get $$ F(a \circ b ) =\frac12 \Big(F(a) b + F(b) a + a \tau(b) + b \tau (a) \Big), \hbox{ for all } a,b\in \mathcal{A}.$$ In this definition the Jordan structure is only considered on the domain algebra. It is worthwhile to mention that if $\mathcal{A}$ and $\mathcal{M}$ are both unital, by taking $b= \mathbf{1}$ in the previous identity we get $F(a) = F(\mathbf{1}) a +\tau(a)$ ($a\in \mathcal{A}$), and if $\tau$ is, in fact, a derivation (as in the case of C$^*$-algebras), $F$ is a generalized derivation of the first type. \smallskip

Motivated by the results in the previous section, we shall say that a linear mapping $F:\mathcal{A} \rightarrow \mathcal{M}$ is a \emph{generalized Jordan derivation of the third type} if there exist (non-necessarily linear) maps $G,H :\mathcal{A} \rightarrow \mathcal{M}$ such that $$ F( a\circ b) = G(a) \circ b + a \circ H(b), \hbox{ for all } a,b\in \mathcal{A}.$$

The following theorem establishes the automatic continuity of every generalized Jordan derivation of the third type when the domain is a C$^*$-algebra. % by adapting the proof of Theorem~\ref{t automatic cont of gener der of type3}. 

\begin{theorem}\label{t automatic cont generalize Jordan der of type 3} Let $\mathcal{A}$ be a C$^{\ast}$-algebra and let $\mathcal{M}$ be a Banach $\mathcal{A}$-bimodule. Then every generalized Jordan derivation of the third type $F:\mathcal{A} \rightarrow \mathcal{M}$ is continuous.
\end{theorem}

\begin{proof} Let $G,H :\mathcal{A} \rightarrow \mathcal{M}$ be two (non-necessarily linear) maps satisfying $$ F( a\circ b) = G(a) \circ b + a \circ H(b), \hbox{ for all } a,b\in \mathcal{A}.$$ Since the Jordan product is commutative, up to replacing $G$ and $H$ by $\frac12 (G+H)$ we can assume that $G=H$. %We have \begin{equation}\label{eq F double G}  F( a\circ b) = G(a) \circ b + a \circ G(b), \hbox{ for all } a,b\in \mathcal{A}.
%\end{equation}
By a new application of Cuntz' theorem \cite[Theorem 1.1 or Corollary 1.2]{Cuntz1976}, we may also assume in the argument concerning the continuity of $F$ that $\mathcal{A}$ is commutative.\smallskip
	 	
Set $$I = \{a \in \mathcal{A} \ : \hbox{ the mapping } x\mapsto F ( a x ) = F (a\circ x) \hbox{ is continuous}\},$$ $$ \hbox{ and }\ J = \{a \in \mathcal{A} \ :  \hbox{ the mapping } x\mapsto a\circ G (x)  \hbox{ is continuous}\}.$$	The identity $F (a x ) = F(a \circ x) = G(a) \circ x + a \circ G(x)$, shows that $I = J$. If $a\in I$ and $b\in \mathfrak{A}$, the mapping $x\mapsto  F ((a b) x) = F ((b a) x) = F (b (a x)) = F (a (bx))$ is clearly continuous. Therefore $ a b = ba \in I$, which shows that $I$ is an ideal of $\mathcal{A}$.\smallskip

To show that $I$ is norm-closed, we take $(a_n)_n\subset I$ converging to some $a\in \mathcal{A}$ in norm. Since by assumptions $a_n \circ G(\cdot)$ is a bounded linear operator from $\mathcal{A}$ to $\mathcal{M}$ which converges pointwise to the linear mapping $a \circ G(\cdot)$, the latter must be bounded by the uniform boundedness principle, and hence $a\in I$.\smallskip

We are now in a position to apply Lemma~\ref{l separate cont} to guarantee that $F|_{I} : {I}\to  \mathcal{M}$ is continuous. The continuity of $F$ will follow from Lemma~\ref{l continuity on a closed finite-codimensional subspace} if we show that ${\mathcal{A}}/{I}$ is finite-dimensional. Otherwise, as in the proof of Theorem~\ref{t automatic cont of gener der of type3}, there exists an infinite sequence $(b_n +I)_n$ of mutually orthogonal non-zero positive elements in ${\mathcal{A}}/{I}$, that is, $(b_n + I) (b_k+ I) = 0,$ for all $n \neq k$, $b_n\geq 0$ and $b_n\neq 0$ for all $n$  (see \cite[Exercise 4.6.13]{K-R}). Choose, a sequence $(c_n)_n$ of mutually orthogonal norm-one positive elements in  $\mathcal{A}$ satisfying $c_n+I = b_n+I$, $\|  c_n^2+I \| = \|  (b_n+I)^2 \| = \| b_n+I\|^2 \neq 0$, for all $n$ (cf. \cite[Exercise 4.6.20]{K-R}). Clearly, the element $c_n^2$ does not belong to $I$, and hence the mapping $x\mapsto F(c_n^2 x)$ is unbounded. We can therefore pick $d_n\in \mathcal{A}$ satisfying $\|d_n\| \leq 2^{-n},$ and $\| F ({{c^2_n}} d_n)\| > \|G(c_n)\| + n$ for all natural $n$. By assumptions, the limit  $a_0 = \sum_{n \geq 1} c_n d_n$ belongs to $\mathcal{A},$ $\|a_0\|\leq 1,$ and $ a_0  c_{m}= c_{m} a_0 = c_m^2 d_m$ for all $m\in \mathbb{N}$. Our hypotheses give  
\begin{align*}
	\infty> 2 \|G(a_0)\| &\geq  \| c_m \circ G (a_0)\|  = \|F(c_m a_0) - G(c_m)\circ  a_0\|  \geq \|F(c_m a_0)\| - \|G(c_m) \circ a_0\| \\ & = \|F(c_m^2 d_m)\| - \|G(c_m) a_0\| \geq m + \|G(c_m)\| - \|G(c_m)\| = m,
\end{align*} for all natural $m$, which is impossible.
\end{proof}

Under some extra hypotheses on the bimodule, we can conclude that every generalized Jordan derivation of the third type on a C$^*$-algebra is, in fact, a generalized derivation of the first, second, and third type.  

\begin{proposition}\label{central elements generalized Jordan der} Let $\mathcal{A}$ be a C$^{\ast}$-algebra and let $\mathcal{M}$ be an essential Banach $\mathcal{A}$-bimodule. Then every generalized Jordan derivation of the third type $F:\mathcal{A} \rightarrow \mathcal{M}$ is a generalized derivation of the second type (and, of course, of the first and third type), that is, there exist $\xi \in \mathcal{M}^{**}$ such that
	$$ F(ab) = F(a)b + aF(b) -a \xi b \ \ \  (\forall a, b \in  \mathcal{A} ). $$
	Furthermore, if $\mathcal{M} $ coincides with $\mathcal{A}$ or $\mathcal{A}^{*}$, the element $\xi$ can be chosen satisfying $\xi a = a \xi$, for all $a\in \mathcal{A}$.
\end{proposition}
  
\begin{proof} We can assume, without loss of generality, that there is a (non-necessarily linear) mapping $G: \mathcal{A}\to \mathcal{M}$ such that \begin{equation}\label{eq f and G} F(a \circ b)  = G(a)\circ b + a \circ G(b), 
	\end{equation} for all $a,b\in \mathcal{A}$ (compare the proof of Theorem~\ref{t automatic cont generalize Jordan der of type 3}).\smallskip  
	
Since $\mathcal{M}$ is assumed to be essential, we can prove that $F$ is a generalized derivation of the second type by combining that $F$ is continuous (cf. Theorem~\ref{t automatic cont generalize Jordan der of type 3}) and Theorem~\ref{t characterization bli generalized derivations}. To see this, just observe that given $a,b\in \mathcal{A}_{sa}$ with $ab = 0$, we can take $d\in \mathcal{A}$ with $d^2 = a $ and $d  b = b d =0$. Therefore, $$b F(a) b= b F(d^2) b= b\left( 2 d\circ G(d)  \right) b =0 ,$$ which concludes the proof of the first statement. %Thus, there exists $\xi\in \mathcal{M}^{**}$ such that $$F( a b ) = F(a) b + a F(b) - a \xi b, \hbox{ for all } a,b\in \mathcal{A},$$ which polarized gives   \begin{equation}\label{eq gener derivation fla on Oct 10 2024} F( a \circ b ) = F(a) \circ b + a \circ F(b) - \frac{a \xi b + b\xi a}{2}, \hbox{ for all } a,b\in \mathcal{A}. \end{equation} Observe that replacing $\xi$ with $\mathbf{1} \xi \mathbf{1}$ we can clearly assume that $\xi= \mathbf{1} \xi \mathbf{1}.$
\smallskip

We shall next show that $G$ is actually linear and continuous. Namely, we can apply similar ideas to those in the proof of Proposition~\ref{p continuous gd 3type}, to see that the mappings of the form $x\mapsto (G(\cdot)\circ b)(x) = G(x)\circ b = F(x\circ b) - x\circ G(b)$ are linear and continuous on $\mathcal{A}$. Since for every $a,b\in \mathcal{A}$ with $\|b\|\leq 1$ we have $\| G(a) \circ b\| = \|F(a \circ b) - a\circ G(b)\| \leq \|F\| \|a\| + \|G(\cdot) \circ a\|,$ the family $\{ G(\cdot) \circ b : \|b\|\leq 1 \}$ is pointwise bounded, and hence uniformly bounded by the uniform boundedness principle. Thus, there exists $K>0$ such that $\| G(\cdot) \circ b \|\leq K$ for all $\|b\|\leq 1$. Let $(u_j)_j$ be a bounded approximate unit in $\mathcal{A}$. For each $a\in \mathcal{A}$, the net $F(a \circ u_j)$ converges to $F(a)$ in norm by the continuity of $F$, and the net $G(a)\circ u_j$ tends to $G(a)$ in norm. Consequently $G$ is linear. We also know that $\| G(a) \| = \lim_{j} \| G(a) \circ u_j\| \leq K \|a\|,$ note that the first equality holds because $\mathcal{M}$ is essential (see, for example, \cite[Theorem 32.22]{HewRosVolII} and \cite[Theorem in page 108, \S~1.15]{CLM1979}). Thus, $G$ is continuous.\smallskip

Having in mind that $F$ is a continuous generalized derivation of the second type on $\mathcal{A}$, Theorem 2.11 in \cite{AyuKudPe2014} assures that $F^{**}$ is a generalized derivation of the second type on $\mathcal{A}^{**},$ and thus there exists ${\xi}\in \mathcal{M}^{**}$ satisfying ${\xi}= \mathbf{1} {\xi} \mathbf{1},$ and \begin{equation}\label{eq F is a generalized derivation of the third type last proposition} F^{**}( a \circ b) = F^{**}( a) \circ b + a\circ F^{**}( b ) - \frac12 (a {\xi} b + b {\xi} a),\hbox{ for all } a,b\in \mathcal{A}^{**}. 
\end{equation} Observe that, a priori, the element ${\xi}$ belongs to $\mathcal{M}^{****}$, however, since $\mathbf{1} F^{**} (\mathbf{1}) \mathbf{1} = \mathbf{1} \widehat{\xi} \mathbf{1} =\widehat{\xi}$ lies in $\mathcal{M}^{**}$.\smallskip

Let us prove the final statement. Assume first that $\mathcal{M}= \mathcal{A}$. The weak$^*$-density of $\mathcal{A}$ in $\mathcal{A}^{**}$ and the separate weak$^*$-continuity of the product in $\mathcal{A}^{**}$ combined with the identity in \eqref{eq f and G} prove that \begin{equation}\label{eq F in terms of G last proposition 1 } F^{**} (a \circ b)  = G^{**} (a)\circ b + a \circ G^{**}(b), \hbox{ for all  } a,b\in \mathcal{A}^{**}.
\end{equation}

Assume next that $\mathcal{M} = \mathcal{A}^{*}$. A consequence of the Grothendieck's inequality, established by Haagerup in \cite[page 95]{Haagerup1985}, asserts that every bounded linear mapping from a C$^*$-algebra to the dual space of another C$^*$-algebra factors through a Hilbert space, and thus it is weakly compact.  Therefore, by Gantmacher's theorem \cite{Gantmacher}, $F^{**}\left(\mathcal{A}^{**}\right)\subseteq \mathcal{M}^{*}$ and  $G^{**}\left(\mathcal{A}^{**}\right)\subseteq \mathcal{M}^{*}$. We claim that \begin{equation}\label{eq F in terms of G in last proposition} F^{**} (a \circ b)  = G^{**} (a)\circ b + a \circ G^{**}(b)
\end{equation} for all $a,b\in \mathcal{A}^{**}$. We just need to extend \eqref{eq f and G} to $\mathcal{A}^{**}$. The left hand side can be treated by a standard argument, since  $F^{**} (a \circ b)$ can be approached by a double weak$^*$-limit via the weak$^*$-continuity of $F^{**}$ and the separate weak$^*$-continuity of the product of $\mathcal{A}^{**}$. Given $b\in \mathcal{A}^{**}$, we can find a net $(b_j)_j\subseteq \mathcal{A}$ such that $G(b_j)\to G^{**} (b)\in\mathcal{A}^{*}$ in the weak$^*$-topology of $\mathcal{A}^{***}$. Thanks to this, we shall handle the summands on the right hand side with the next property: if $(a_j)_j$ and $(\phi_i)$ are two nets in $\mathcal{A}$ and $\mathcal{A}^*$ converging to $a\in \mathcal{A}^{**}$ and $\phi\in \mathcal{A}^*$ in the $\sigma(\mathcal{A}^{**},\mathcal{A}^{*})$ and the  $\sigma(\mathcal{A}^{***},\mathcal{A}^{**})$ topologies,  respectively, we have $$\frac12 \lim_{j} \lim_{i} (a_j  \phi_i +  \phi_i a_j ) = \lim_{j} \lim_{i} (a_j \circ \phi_i) = a\circ \phi = \frac12 (a \phi + \phi a) \hbox{ in the $\sigma(\mathcal{A}^{***},\mathcal{A}^{**})$ topology.}$$ The basic properties of the module operations in $\mathcal{A}^{***}$ assure that $\lim_{j} \lim_{i} a_j  \phi_i =a \phi $ in the $\sigma(\mathcal{A}^{***},\mathcal{A}^{**})$ topology, and $w^*\hbox{-}\lim_{i}  \phi_i a_j=  \phi a_j$ for all $j$.  Fix now $c\in \mathcal{A}^{**}$. The net $(a_j c)_j$ converges to $a c$ in the weak$^*$-topology of $\mathcal{A}^{**}$. By applying that $\phi\in \mathcal{A}^{*}$, we deduce that $$\lim_{j} (\phi a_j ) (c) = \lim_{j} \phi (a_j  c) = \phi (a c) = (\phi a) (c),$$ and thus $w^*\hbox{-}\lim_{j}  \phi a_j  = \phi a$, which finishes the proof of \eqref{eq F in terms of G in last proposition}.\smallskip

Summarizing, by \eqref{eq F is a generalized derivation of the third type last proposition}, \eqref{eq F in terms of G last proposition 1 }, and \eqref{eq F in terms of G in last proposition} we arrive to $G^{**} (a) = F^{**}(a) - a\circ G^{**} (\mathbf{1})$, $ G^{**} (\mathbf{1}) = \frac12 F^{**} (\mathbf{1}) = \frac12 \xi$, 
$$\begin{aligned}
	& F^{**}( a) \circ b + a\circ F^{**}( b ) - \frac12 (a \xi b + b \xi a) = F^{**} (a\circ b ) \\
	&= F^{**} (a) \circ b + a\circ F^{**} ( b ) - (a\circ G^{**} (\mathbf{1})) \circ b - (b\circ G^{**} (\mathbf{1})) \circ a,
\end{aligned}$$ equivalently, $$ a \xi b + b \xi a = (a\circ \xi) \circ b + (b\circ \xi) \circ a \Leftrightarrow  \xi a b + b a\xi -2 b \xi a + \xi  b a + a b \xi -2 a \xi b =0, $$ for all $a,b\in \mathcal{A}^{**}$. In particular, $ \xi p +  p \xi = 2 p \xi p$, for every projection $p\in \mathcal{A}^{**}$, which proves that $\xi p =p \xi p = p \xi$, and hence $\xi$ commutes with all elements in $\mathcal{A}^{**}$.
\end{proof}

\begin{remark}\label{third type is stictly included in second} We can conclude now that for a non-commutative C$^*$-algebra $\mathcal{A}$, the class of generalized Jordan derivations of the third type on $\mathcal{A}$ is strictly included in the class of generalized Jordan derivations of the second type. Consider, for example, an element $\zeta$ which is not in the centre of  $\mathcal{A}$ and the mapping $F: \mathcal{A}\to \mathcal{A}$, $F(a )= a\circ  \zeta$ which is a generalized Jordan derivation of the second type but not of the third type. 
\end{remark}

Finally, it is worth to say a few words about what could be the Jordan version of generalized derivations of the first type. Let $\mathcal{M}$ be a Banach $\mathcal{A}$-bimodule on a Banach algebra $\mathcal{A}$. A \emph{generalized Jordan derivation of the first type} is a linear mapping $F : \mathcal{A}\to \mathcal{M},$ for which there exists a Jordan derivation $D:   \mathcal{A}\to \mathcal{M}$ satisfying $F( a\circ b) = F(a) \circ b + a \circ D(b)$ for all $a,b\in \mathcal{A}$. Indeed some algebraists defined the so-called ``Jordan generalized derivation" by letting $D$ be just a linear map. That is the case in the paper by Li and Benkovi\v{c} \cite{LiB2011}, where it is shown that any generalized Jordan derivation of the first type on a triangular algebra is a kind of generalized derivation of the third type. Every Jordan derivation of the first type is automatically a generalized Jordan derivation of the third type, and thus continuous by Theorem~\ref{t automatic cont generalize Jordan der of type 3} when $\mathcal{A}$ is a C$^*$-algebra. %Furthermore, if $\mathcal{M}$ is essential, the sets of all generalized Jordan derivations of the first and third type are both strictly included in the set of all generalized Jordan derivations of the second type, and actually they are all generalized derivations of the first, second, and third type. 
\smallskip

Despite we have tried to clarify the relationships between the different notions of ``generalized derivations'', several interesting questions remain open after our study.  Can we remove the hypothesis $\mathcal{M}$ being essential in Proposition \ref{central elements generalized Jordan der}? Is the final conclusion in the same proposition true for other Banach bimodules besides $\mathcal{A}$ and $\mathcal{A}^*$? \smallskip

It is also natural to expect that the notions of generalized derivations of the third type studied in this note give rise to appropriate concepts in the setting of JB$^*$-algebras and JB$^*$-triples. In the case of a unital JB$^*$-algebra $\mathcal{J}$,% and a Jordan Banach module $\mathcal{M}$,
 the techniques and arguments in \cite{HejNik96} can be appropriately modified to prove that every generalized Jordan derivation of the third type on $\mathcal{J}$ or from $\mathcal{J}$ to $\mathcal{J}^*$ is automatically continuous. However, under weaker hypotheses the conclusion is not so clear. \medskip

%%%%%%%%%%%%%%%%%%%%%%%%%%
%\vspace{1cm}
\textbf{Acknowledgements} Second and third authors partially supported by grant \linebreak PID2021-122126NB-C31 funded by MICIU/AEI/10.13039/501100011033 and by ERDF/EU. Second author also supported by Junta de Andalucía grant FQM375, the IMAG--Mar{\'i}a de Maeztu grant CEX2020-001105-M/AEI/10.13039/501100011033, and by the (MOST) Ministry of Science and Technology of China grant G2023125007L. Third author supported by China Scholarship Council Program (Grant No.202306740016).\smallskip

\textbf{Disclosure:} All authors declare that they have no conflicts of interest to disclose.

%-----------------------------------------------------------------------------
%-----------------------------------------------------------------------------
\end{document}